\documentclass[12pt]{article}
\usepackage{amsthm,amsfonts,amssymb,amscd}

\begin{document}

%\begin{flushleft}
%ÓÄÊ 512.7
%\end{flushleft}

\begin{center}
\large \bf  $K$-trivial structures \\ on Fano complete
intersections
\end{center}\vspace{0.5cm}

\begin{center}
Aleksandr Pukhlikov
\end{center}\vspace{0.5cm}

\parshape=1
3cm 10cm \noindent {\small \quad\quad\quad \quad\quad\quad\quad
\quad\quad\quad {\bf }\newline It is proven that any structure of
a fiber space into varieties of Kodaira dimension zero on a
generic Fano complete intersec\-tion of index 1 and dimension $M$
in ${\mathbb P}^{M+k}$ for $M\geq 2k+1$ is a pencil of hyperplane
sections. We describe $K$-trivial structures on varieties with a
pencil of Fano complete intersec\-tions.

\noindent Bibliography: 9 items.} \vspace{1cm}

\noindent {\bf 1. Formulation of the main result.} Fix $k\geq 2$,
$M\geq 2k+1$ and a set of integers $(d_1,\dots,d_k)\in{\mathbb
Z}^k_+$ satisfying the conditions $d_k\geq\dots\geq d_1\geq 2$ and
$$
d_1+\dots+d_k=M+k.
$$
Assume in addition that if $M=2k+2$, then $k\geq 4$, and if
$M=2k+1$, then $k\geq 5$ and $(d_1,\dots,d_k)\neq
(2,\dots,2,k+3)$.

The symbol ${\mathbb P}$ denotes the complex projective space
${\mathbb P}^{M+k}$. Consider a smooth complete intersection of
codimension $k$ in ${\mathbb P}$
$$
V=\{f_1=\dots=f_k=0\}\subset{\mathbb P},
$$
where $f_i\in H^0({\mathbb P},{\cal O}_{\mathbb P}(d_i))$ are
homogeneous polynomials of degree $d_i$. The aim of this note is
to give a complete description of structures of non-maximal
Kodaira dimension on a generic complete intersection $V$, from
which one immediately derives a description of $K$-trivial
structures on varieties with a pencil of complete intersections.\vspace{0.1cm}

Let $\beta\colon W\to S$ be a morphism of projective varieties
with connected fibers (a fiber space). By the {\it relative
Kodaira dimension} $\kappa(W/S)$ of the fiber space $W/S$ we mean
the Kodaira dimension of a fiber of general position
$\beta^{-1}(s)$, $s\in S$. By a {\it structure} of a fiber space
of the relative Kodaira dimension
$\kappa\in\{-\infty,0,1,\dots,M\}$ on the variety $V$ we mean an
arbitrary birational map $\chi\colon V\dashrightarrow W$, where
$\beta\colon W\to S$ is a fiber space of the relative Kodaira
dimension $\kappa$. It is well known [1], that the generic
complete intersection $V$ is birationally superrigid, in
particular, on $V$ there are no structures of negative relative
Kodaira dimension. The main result of this note is\vspace{0.1cm}

{\bf Theorem 1.} {\it Let $\chi\colon V\dashrightarrow W$ be a
structure of a fiber space of non-maximal relative Kodaira
dimension on a sufficiently general complete intersection
$V\subset{\mathbb P}$, that is, the inequality
$$
\kappa(W/S)+\mathop{\rm dim}S<\mathop{\rm dim}W=M
$$
holds. Then $\kappa(W/S)=0$, $S={\mathbb P}^1$ and there exists a
uniquely determined linear subspace $\Lambda\subset{\mathbb P}$
of codimension two such that the following diagram of maps is
commutative
$$
\begin{array}{ccccc}
& V & \stackrel{\chi}{\dashrightarrow} & W &\\
\pi & \downarrow & & \downarrow & \beta\\
& {\mathbb P}^1 & = & S &,\\ \end{array}
$$
where $\pi=\pi_{\Lambda}|_V\colon V\dashrightarrow{\mathbb P}^1$
is the restriction onto $V$ of the linear projection
$\pi_{\Lambda}\colon{\mathbb P}\dashrightarrow{\mathbb P}^1$ from
the subspace $\Lambda$}.\vspace{0.1cm}

Obviously, the restriction onto $V$ of the projection
$\pi_{\Lambda}$ from a subspace of codimension two is a
$K$-trivial structure. Theorem 1 immediately implies a complete
description of such structures on varieties with a pencil of Fano
complete intersections. Recall the construction of such varieties
[2].\vspace{0.1cm}

Let $\Pi\stackrel{\mu}{\to}{\mathbb P}^1$ be a locally trivial
bundle with the fiber ${\mathbb P}$. Consider a smooth subvariety
$X\subset\Pi$ of codimension $k$, such that every fiber of the
projection $\mu_X=\mu|_X\colon X\to{\mathbb P}^1$ is a (possibly
singular) Fano complete intersection of the type
$d_1\cdot\dots\cdot d_k$ in $\mu^{-1}(t)\cong{\mathbb P}$. The
variety $X$ is assumed to be generic in the sense of regularity
conditions [2,Theorem 5]; besides, the fiber of general position
$\mu^{-1}_X(t)$ is assumed to be generic in the sense of
Proposition 1, which is formulated and proved below.\vspace{0.1cm}

Let $\sigma\colon{\mathbb P}^1\to\mathop{\rm Gr}(M+k-2,\Pi)$ be
an arbitrary map, associating to a point $t\in{\mathbb P}^1$ a
linear subspace of codimension two $\sigma(t)\in\mathop{\rm
Gr}(M+k-2,\mu^{-1}(t))$ in the fiber $\mu^{-1}(t)\cong{\mathbb
P}$. By the symbol $\Delta_{\sigma}$ we denote the ruled surface,
consisting of the hyperplanes $H\subset\mu^{-1}(t)$, containing
the subspace $\sigma(t)$, by the symbol
$$
\pi_{\sigma}\colon\Pi\dashrightarrow\Delta_{\sigma}
$$
the fiber-wise projection, $\pi_{\sigma}|_{\mu^{-1}(t)}$ is the
projection from the subspace $\sigma(t)$ onto ${\mathbb P}^1$.
Theorem 1 implies directly\vspace{0.1cm}

{\bf Theorem 2.} {\it Assume that a Fano fiber space $X/{\mathbb
P}^1$ satisfies the conditions

{\rm (i)} $K^2_X+2H_F\not\in\mathop{\rm Int}A^2_+X$,

{\rm (ii)} $-K_X\not\in A^1_{\rm mov}X$,

\noindent where $H_F=(-K_X\cdot F)$ is the class of a hyperplane
section of the fiber of the projection $\mu_X$ in $A^2X$. Then
for any structure $\chi\colon X\dashrightarrow W$ of a fiber
space of the relative Kodaira dimension zero we get $\mathop{\rm
dim}S=2$ and, moreover, there are a section $\sigma\colon{\mathbb
P}^1\to\mathop{\rm Gr}(M+k-2,\Pi)$ and a birational map
$\gamma\colon\Delta_{\sigma}\dashrightarrow S$ such that the
following diagram is commutative:}
$$
\begin{array}{ccccc}
& X & \stackrel{\chi}{\dashrightarrow} & W &\\
\pi_{\sigma} & \downarrow & & \downarrow & \beta\\
& \Delta_{\sigma} & \stackrel{\gamma}{\dashrightarrow} & S. &\\
\end{array}
$$

\noindent One obtains Theorem 2 from Theorem 1 in a trivial way, see the
proof of Theorem 5 in [2].\vspace{0.3cm}

\newpage

%%%%%%%%%%%%%%%%%%%%%%%%%%%%%%%%%%%%%%%%%%%%%%%%%%%%%%%%%%%%%%%%%%%
%%%%%%%%%%%%%%%%%%%%%%%%%   section 2

\noindent {\bf 2. A summary of known results.} It was observed in
XIX century that the problems of description of structures of a
rationally connected fiber space (or structures of negative
relative Kodaira dimension) and of structures of relative Kodaira
dimension zero are parallel to each other. If for a given
rationally connected variety the first problem admits a complete
solution, then the second problem can be solved by the same
methods (but not the other way round!). In the modern period of
algebraic geometry the first paper describing the structures of
relative Kodaira dimension zero was [3]. In [4] a description of
$K$-trivial structures on generic Fano hypersurfaces of index 1
(an analog of Theorem 1 for these varieties) was derived from the
results of [5]. See also [6,7]. In [1] a sketch of the proof of
the following fact was given.\vspace{0.1cm}

{\bf Theorem 3.}  {\it In the assumptions of Theorem 1 every
structure $\chi\colon V\dashrightarrow W$ of relative Kodaira
dimension zero on $V$ is a pencil:} $\mathop{\rm dim}S=1$.\vspace{0.1cm}

Theorem 3 is weaker than Theorem 1. In fact, in [1, Sec. 0.4] a
proof of the following fact was sketched: in the assumptions of
Theorem 1, let $\Sigma$ be the $(\beta\circ \chi)$-preimage on
$V$ of an arbitrary movable linear system on the base $S$,
$\Sigma\subset|-nK_V|$ for some $n\geq 1$. Then there exists an
irreducible subvariety of codimension two $B\subset V$ such that
\begin{equation} \label{1}
\mathop{\rm mult}\nolimits_B\Sigma=n,
\end{equation}
in particular, the self-intersection of this movable linear system
$Z=(D_1\circ D_2)$, $D_i\in\Sigma$, is $Z=n^2B$ (whence Theorem 3
immediately follows). Since the subvariety $B$ is numerically
equivalent to a section of the variety $V$ by a linear subspace
of codimension two, in [1, Sec. 0.4] it was conjectured that
$\Sigma$ is composed from some pencil of hyperplane sections of
$V$ (this is equivalent to the claim of Theorem 1). In [8] a more
detailed proof of Theorem 3 was given, using the cone techniques
[5,9]. However, in the argument, given in [8], there was a gap in the proof of the lemma (the only lemma in that paper), generalizing the cone technique. (Namely, to describe the intersection of the base of a cone and a curve lying on the cone, a set of free linear systems on the cone was considered. Taking a generic divisor in each system, it was concluded that the zero-dimensional intersection of all those divisors with the base of the cone consists of distinct points, that is, each point in the intersection is of multiplicity one. However, such a conclusion is not obvious and requires a proof. Usually such claims are deduced from the Bertini theorem, but in the situation under consideration the Bertini theorem can not be used without a special justification. The Bertini theorem claims that the singularities of a generic
divisor of a movable linear system are concentrated on its base
set. In [8] a set of movable
linear systems is considered, however, the divisors in those systems are chosen
not independently of each other, more precisely, once a generic
divisor in one of them has been chosen, the choice of the
remaining divisors is uniquely determined. Therefore, the Bertini
theorem does not apply or, in order to apply it, one needs a
special argument.) Finally, in [2] the argument of [8] was
replaced by another one and the proof of Theorem 3 was completed
(see [2, Proposition 3.5]). As an immediate corollary, the
following fact was obtained.\vspace{0.1cm}

{\bf Theorem 4.} {\it In the assumptions of Theorem 2 for any
structure $\chi\colon X\dashrightarrow W$ of relative Kodaira
dimension zero we get $\mathop{\rm dim}S=2$ and the structure
$\chi$ is compatible with the projection} $X\to{\mathbb P}^1$.
(See [2, Theorem 5].)\vspace{0.1cm}

Theorem 4 is a weaker version of Theorem 2. In this note, we make
the concluding step in the description of structures of zero
Kodaira dimension on Fano complete intersections.\vspace{0.3cm}

%%%%%%%%%%%%%%%%%%%%%%%%%%%%%%%%%%%%%%%%%%%%%%%%%%%%%%%%%%%%%%%%%%%
%%%%%%%%%%%%%%%%%%%%%%%%%   section 3

\noindent {\bf 3. The structure of the proof of Theorem 1.} The
variety $V$ is assumed to be generic in the following sense.
Firstly, $V$ satisfies the regularity condition [1] at every
point, see [1, Sec. 1.2]. Secondly, the following claim holds.\vspace{0.1cm}

{\bf Proposition 1.} {\it A sufficiently general complete
intersection $V$ satisfies the following condition: for any
linear subspace $\Lambda\subset{\mathbb P}$ of codimension $k+1$
the intersection $V\cap\Lambda$ is an irreducible reduced variety
of dimension} $M-k-1$.\vspace{0.1cm}

{\bf Proof} is given in Sec. 5. By the Lefschetz theorem the claim
of the proposition holds for $M\geq 2k+3$ for any smooth complete
intersection $V$. Only two cases need to be considered: $M=2k+2$
and $M=2k+1$.\vspace{0.1cm}

We assume that the complete intersection $V$ satisfies the
property of Proposition 1.\vspace{0.1cm}

Let $\chi\colon V\dashrightarrow W$ be a structure of a fiber
space of non-maximal relative Kodaira dimension, $\Sigma$ the
strict transform on $V$ with respect to $\chi$ of the pull back
on $W$ of some movable linear system on the base $S$. Then
$\Sigma\subset|-nK_V|$, $n\geq 1$, the pair
$(V,\frac{1}{n}\Sigma)$ is not terminal, the system $\Sigma$ is
composed from a pencil and there is an irreducible subvariety
$B\subset V$ of codimension two, such that the inequality
(\ref{1}) holds [1, Proposition 3.5]. Let $b\in B$ be a point of
general position, $\Delta=T_bB\subset{\mathbb P}$ the tangent
space.\vspace{0.1cm}

{\bf Proposition 2.} {\it Let $\mu_{\Delta}\colon{\mathbb
P}\dashrightarrow{\mathbb P}^{k+1}$ be the linear projection from
$\Delta$,
$$
\mu=\mu_{\Delta}|_V\colon V\dashrightarrow{\mathbb
P}^{k+1}
$$
its restriction onto $V$. Then the fibers of the rational map
$\mu$ are irreducible and reduced, whereas the linear system
$\Sigma$ is the pull back via $\mu$ of a movable linear system
$\Gamma$ on ${\mathbb P}^{k+1}$}.\vspace{0.1cm}

{\bf Proof} is given in Sec. 4.\vspace{0.1cm}

{\bf Proof of Theorem 1.} By Proposition 1, the set $V\cap\Delta$
has codimension $k+2$ in $V$. Considering the dimensions, we
conclude that $\mu(B)=\bar B$ is an irreducible subvariety of
codimension two ($\bar B$ can not be a divisor, since in that
case the linear system $\Gamma$, and therefore also $\Sigma$,
would have had a fixed component), so that $B=\mu^{-1}(\bar B)$.
Therefore,
$$
\mathop{\rm mult}\nolimits_{\bar B}\Gamma=n,
$$
where $\Gamma$ is a linear system of hypersurfaces of degree $n$.
This is possible in one case only, when $\bar B\subset{\mathbb
P}^{k+1}$ is a linear subspace of codimension two, and the system
$\Gamma$ is composed from the pencil of hyperplanes, containing
$\bar B$. But then $B=\Lambda\cap V$, where
$\Lambda=\mu^{-1}_{\Delta}(\bar B)$ is a linear subspace of
codimension two in ${\mathbb P}$, and the system $\Sigma$ is
composed from the pencil of sections of $V$ by hyperplanes
containing $\Lambda$. Q.E.D. for Theorem 1.\vspace{0.1cm}

As we noted above, Theorem 2 follows immediately from Theorem
1.\vspace{0.3cm}

%%%%%%%%%%%%%%%%%%%%%%%%%%%%%%%%%%%%%%%%%%%%%%%%%%%%%%%%%%%%%%%%%%%
%%%%%%%%%%%%%%%%%%%%%%%%%   section 4

\noindent {\bf 4. Linear projections and cones.} Let us prove
Proposition 2. Proposition 1 implies that the fibers of $\mu$ are
irreducible and reduced. To prove the main claim that
$\Sigma=\mu^*\Gamma$, let us consider a point of general position
$p\in V$. Let $D\in\Sigma$ be the divisor, containing that point.
Now we get\vspace{0.1cm}

{\bf Proposition 3.} {\it The following inclusion holds}
$$
T_p\mu^{-1}(\mu(p))\subset T_pD.
$$

Proposition 2 follows immediately from Proposition 3: as the
point $p$ is generic, we get that $\mu(D)\subset{\mathbb P}^{k+1}$
is a divisor (that is, $\mu(D)\neq{\mathbb P}^{k+1}$), so that
$\Sigma=\mu^*\Gamma$ for some linear system $\Gamma$, which is
what we need. Q.E.D.\vspace{0.1cm}

{\bf Proof of Proposition 3.} Set $P=\langle\Delta, p\rangle$ to
be the fiber of the linear projection $\mu_{\Delta}$, so that
$\mu^{-1}(\mu(p))=P\cap V$. Since the point $p$ is generic and
the fiber $\mu^{-1}(\mu(p))$ is non-singular at this point, we get
$$
\mathop{\rm codim}\nolimits_{\mathbb P}(P\cap T_pV)= \mathop{\rm
codim}\nolimits_{\mathbb P}P+k,
$$
that is, the linear subspaces $P$ and $T_pV\subset{\mathbb P}$
are in general position.\vspace{0.1cm}

Assume now that the point of general position $b\in B$ was chosen
in the following way: on the variety $B$ we considered an
arbitrary family of irreducible $k$-dimensional subvarieties
$\{Y_u,u\in U\}$, sweeping out $B$, in that family we chose a
variety $Y=Y_u$ of general position, and the point $b$ is a point
of general position on $Y$. In particular, $T_bY$ is a generic
linear subspace of dimension $k$ in $\Delta=T_bB\subset {\mathbb
P}$. In particular,
\begin{equation}\label{2}
\dim(\langle T_bY,p\rangle\cap T_pV)=1,
\end{equation}
that is, the linear subspaces $\langle T_bY,p\rangle$ è $T_pV$
are in general position.\vspace{0.1cm}

By the symbol $[b,p]$ we denote the line in ${\mathbb P}$,
joining these two points, by the symbol $(b,p)$ the set
$[b,p]\setminus \{b,p\}$. Take a point $x\in (b,p)$. Set $C(Y,x)$
to be the cone with the vertex $x$ and the base $Y$.\vspace{0.1cm}

{\bf Proposition 4.} {\it For sufficiently general $Y$, $b$, $p$,
$x$ the following claims are true:

\noindent {\rm (i)} the point $z\in C(Y,x)$ is a singularity of
that cone, if either $z=x$, or $z\in [y,x]$, where $y\in
\mathop{\rm Sing} Y$,

\noindent {\rm (ii)} the closed algebraic set $R(Y,x)$, which is
the union of all irreducible components of the intersection
$C(Y,x)\cap V$, containing the point $p$, is an irreducible
curve, non-singular at the point $p$,

\noindent {\rm (iii)} the curve $R(Y,x)$ intersects the
subvariety $Y$ outside the closed subset $\mathop{\rm Sing} Y$ of
singular points of this variety.}\vspace{0.1cm}

{\bf Proof} of the claim (i) is given in [2, Sec. 3.3.1] (in
addition to the arguments, given in [2], one needs to note that
$V$ can not be contained in the variety of secant lines
$\mathop{\rm Sec} Y$ of the variety $Y$: even if $M=2k+1$, the
variety $V$ is not covered by lines). Furthermore, $[b,p]\cap
Y=\{b\}$, so that the point $p\in C(Y,x)$ is non-singular.
Obviously,
$$
T_pC(Y,x)=\langle T_bY, p\rangle,
$$
whence, taking into account (\ref{2}), it follows that the
varieties $C(Y,x)$ and $V$ intersect transversally at the point
$p$, which proves (ii). Finally, the claim (iii) is proved by the
arguments of [2, Sec. 3.3.1], taking into account that when $Y$,
$b\in Y$ and $p\in V$ vary, the points $x\in (b,p)$ fill out an
open subset of the projective space $\mathbb P$ (for instance, by
the surjectivity of the map $\mu$). Q.E.D. for Proposition 4.\vspace{0.1cm}

Now we argue as in [2, Sec. 3.3]: the curve $R(Y,x)$ meets $Y$ at
the points, which are non-singular both on $Y$ and on the cone
$C(Y,x)$. Therefore, we get that the intersection number
$(R(Y,x)\cdot Y)_{C(Y,x)}$ is well defined.\vspace{0.1cm}

{\bf Lemma 1.} {\it The following equality holds:} $(R(Y,x)\cdot
Y)_{C(Y,x)}=\mathop{\rm deg} R(Y,x)$.\vspace{0.1cm}

{\bf Proof} is given in [9, \S 1].\vspace{0.1cm}

Now let us come back to the divisor $D\in\Sigma$, containing the
point $p$.\vspace{0.1cm}

{\bf Lemma 2.} {\it The following inclusion holds:} $R(Y,x)\subset
D$.\vspace{0.1cm}

{\bf Proof.} Obviously, $(D\cdot R(Y,x))=n\mathop{\rm deg}
R(Y,x)$. On the other hand, as it was shown in [2, Sec. 3.3.2],
$$
\sum_{y\in D\cap R(Y,x)\cap Y} (D\cdot R(Y,x))_y\geq (R(Y,x)\cdot
Y)_{C(Y,x)}\mathop{\rm mult}\nolimits_Y D.
$$
Taking into account that $\mathop{\rm mult}\nolimits_Y
D=\mathop{\rm mult}\nolimits_B D=n$ and the equality of Lemma 1,
we obtain the claim of Lemma 2, since
$$
p\in D\cap R(Y,x)
$$
and $p\not\in Y$. Q.E.D. for the lemma.\vspace{0.1cm}

Lemma 2 implies the inclusion
$$
\langle T_bY,p\rangle\cap T_pV\subset T_pD,
$$
whence, since the subspace $T_bY\subset \Delta$ is generic, we
get the inclusion
$$
\langle \Delta,p\rangle\cap T_pV\subset T_pD,
$$
which completes the proof of Proposition 3.\vspace{0.3cm}

%%%%%%%%%%%%%%%%%%%%%%%%%%%%%%%%%%%%%%%%%%%%%%%%%%%%%%%%%%%%%%%%%%%
%%%%%%%%%%%%%%%%%%%%%%%%%   section 5

{\bf 5. Complete intersections of general position.} Let us prove
Proposition 1. As we noted in Sec. 3, we have to consider the two
cases: when $M=2k+2$ and $M=2k+1$. Consider first the following
general problem. Let $X\subset {\mathbb P}^{N}$ be an irreducible
subvariety of dimension $l\geq 2$. By the symbol ${\cal
P}_{d}={\cal P}_{d,N}$ we denote the space of homogeneous
polynomials of degree $d$ on ${\mathbb P}^{N}$. Let
$U_d(X)\subset {\cal P}_{d}$ be the open set, consisting of such
polynomials $f\in {\cal P}_{d}$, that $\{f|_X=0\}$ is an
irreducible reduced subvariety of dimension $(l-1)$.
Respectively, let $R_d(X)={\cal P}_{d}\setminus U_d(X)$ be the
set of ``incorrect'' polynomials. The problem is to estimate from
below the codimension of the closed set $R_d(X)$ in the space
${\cal P}_{d}$.\vspace{0.1cm}

{\bf Lemma 3.} {\it The following estimate holds:}
$$
\mathop{\rm codim}(R_d(X)\subset {\cal P}_{d})\geq {d+l-2 \choose
d} -l+1.
$$

{\bf Proof.} Let $\gamma\colon {\mathbb P}^{N}\dashrightarrow
{\mathbb P}^{l-1}$ be the linear projection from a $(N-l)$-plane
of general position, $\gamma_X=\gamma|_X\colon X\dashrightarrow
{\mathbb P}^{l-1}$ its restriction onto $X$ (the set of points
where $\gamma_X$ is not defined is zero-dimensional). Obviously,
the map $\gamma_X$ is surjective, all its fibers are
one-dimensional and the fiber $\gamma_X^{-1}(z)$ over a point of
general position $z\in {\mathbb P}^{l-1}$ is an irreducible
curve. The set $\Delta\subset {\mathbb P}^{l-1}$, consisting of
such points $z$, that the fiber $\gamma_X^{-1}(z)$ is reducible
or non-reduced, is a proper closed subset of ${\mathbb P}^{l-1}$
(at most a divisor). Therefore, for any irreducible divisor $D$ on
${\mathbb P}^{l-1}$, such that $D\not\subset\Delta$, its inverse
image
$$
\gamma_X^{-1}(D)=\gamma^{-1}(D)\cap X
$$
is irreducible and reduced. In other words, for any irreducible
polynomial $f$ on ${\mathbb P}^{l-1}$, such that
$\{f=0\}\not\subset\Delta$, we get $\gamma^*f\in U_d(X)$. The set
$\gamma^*{\cal P}_{d,l-1}$ is a linear subspace of the space
${\cal P}_{d}$ (the same polynomials considered as polynomials in
a larger number of variables). Let $R_{d,l-1}\subset {\cal
P}_{d,l-1}$ be the closed subset of reducible polynomials. From
what was said, it follows that
$$
\mathop{\rm codim} (R_d(X)\subset {\cal P}_{d})\geq \mathop{\rm
codim} (R_{d,l-1}\subset {\cal P}_{d,l-1}).
$$
The last codimension is easy to compute: the irreducible
component of maximal dimension of the set $R_{d,l-1}$ consists of
the polynomials of the form $f=f^{\sharp}h$, where $h\in {\cal
P}_{1,l-1}$ is a linear form. Q.E.D. for Lemma 3.\vspace{0.1cm}

Let us come back to the proof of Proposition 1. Assume that
$M=2k+2$. Assume also that the complete intersection $V$ is
generic in the following sense: the variety
\begin{equation}\label{3}
V^{\sharp}=\{f_1=\dots =f_{k-1}=0\}\subset{\mathbb P}
\end{equation}
is smooth. Obviously, $\dim V^{\sharp}=2k+3$, so that by the
Lefschetz theorem the intersection $V^{\sharp}\cap \Lambda$ is
irreducible, reduced and has dimension $k+2$ for any linear
subspace $\Lambda\subset{\mathbb P}$ of codimension $k+1$. Fix
such a subspace. Obviously,
$$
V\cap \Lambda=\{ f_k|_{V^{\sharp}\cap\Lambda}=0\}.
$$
Let $R_{\Lambda}\subset {\cal P}_{d_k,M+k}$ be the closed subset
of polynomials $f_k$ of degree $d_k$, for which $V\cap\Lambda$ is
not an irreducible reduced subvariety of dimension $k+1$. By
Lemma 3,
\begin{equation}\label{4}
\mathop{\rm codim} R_{\Lambda}\geq {d_k+k\choose d_k} -k-1.
\end{equation}
The equality $d_1+\dots+d_k=M+k$ implies that $d_k\geq 4$.
Elementary computations show that the right hand part of the
inequality (\ref{4}) is strictly higher than the dimension of the
projective Grassmanian of $(2k+1)$-planes in $\mathbb P$ for
$k\geq 4$. This proves Proposition 1 for $M=2k+2$.\vspace{0.1cm}

Assume now that $M=2k+1$. In this case the arguments are
completely similar to those given above, but we need two steps.
First, we consider the complete intersection
$$
V^{+}=\{f_1=\dots =f_{k-2}=0\}\subset{\mathbb P},
$$
which is assumed to be smooth. By the Lefschetz theorem the
intersection $V^+\cap\Lambda$ is irreducible and reduced. Now the
arguments similar to those given above for $M=2k+2$, show that
for a generic polynomial $f_{k-1}$ of degree $d_{k-1}\geq 3$ the
closed set
$$
\{ f_{k-1}|_{V^{+}\cap\Lambda}=0\}
$$
is irreducible and reduced for any subspace $\Lambda$. Now we
consider the variety $V^{\sharp}$, defined by the formula
(\ref{3}), and argue as in the case $M=2k+2$ and complete the
proof. We omit the details of elementary computations. Q.E.D. for
Proposition 1.\vspace{0.1cm}

{\bf Remark 1.} The additional (compared to [1]) restrictions for
the parameters $k$, $d_1$,\dots $d_k$ are needed precisely for
the reason that for the excluded values the proof of Proposition
1 does not work. However, there are no doubts that both the claim
of Proposition 1 and, the more so, Theorem 1 are true for those
values as well. Here is the list of excluded families:\vspace{0.1cm}

$$
2\cdot 5\quad \mbox{and}\quad 3\cdot 4\quad \mbox{in}\quad
{\mathbb P}^{7},
$$
$$
2\cdot 6,\quad 3\cdot 5\quad \mbox{and}\quad 4\cdot 4\quad
\mbox{in}\quad {\mathbb P}^{8},
$$
$$
2\cdot 2\cdot 6,\quad 2\cdot 3\cdot 5,\quad 2\cdot 4\cdot 4\quad
\mbox{and}\quad 3\cdot 3\cdot 4\quad \mbox{in}\quad {\mathbb
P}^{10},
$$
$$
2\cdot 2\cdot 7,\quad 2\cdot 3\cdot 6,\quad 2\cdot 4\cdot 5,\quad
3\cdot 3\cdot 5 \quad\mbox{and}\quad 3\cdot 4\cdot 4\quad
\mbox{in}\quad {\mathbb P}^{11},
$$
$$
2\cdot 2\cdot 3\cdot 6,\quad 2\cdot 2 \cdot 4\cdot 5,\quad 2\cdot
3\cdot 3\cdot 5\quad \mbox{and}\quad 2\cdot 3\cdot 4\cdot 4\quad
\mbox{in}\quad {\mathbb P}^{13},
$$
and the infinite series $2\cdot \dots \cdot 2\cdot (k+3)$ in
${\mathbb P}^{3k+1}$, $k\geq 2$.\vspace{0.5cm}

{\bf References}\vspace{0.3cm}

\noindent 1. Pukhlikov A.V., Birationally rigid Fano complete
intersections, Crelle J. f\" ur die reine und angew. Math. {\bf
541} (2001), 55-79. \vspace{0.3cm}

\noindent 2. Pukhlikov A.V., Birational geometry of algebraic
varieties with a pencil of Fano complete intersections,
Manuscripta Mathematica {\bf 121} (2006), 491-526. \vspace{0.3cm}

\noindent 3. Dolgachev I.V., Rational surfaces with a pencil of
elliptic curves. Izvestiya: Mathematics. {\bf 30} (1966),
1073-1100.\vspace{0.3cm}

\noindent 4. Cheltsov I.A., Log pairs on hypersurfaces of degree
$N$ in ${\mathbb P}^N$. Math. Notes. {\bf 68} (2000), no. 1,
131-138.\vspace{0.3cm}

\noindent 5. Pukhlikov A.V., Birational automorphisms of Fano
hypersurfaces, Invent. Math. {\bf 134} (1998), no. 2, 401-426.
\vspace{0.3cm}

\noindent 6. Cheltsov I.A., Birationally rigid del Pezzo
fibrations. Manuscripta Math. {\bf 116} (2005), no. 4,
385--396.\vspace{0.3cm}

\noindent 7. Cheltsov Ivan and Karzhemanov Ilya, Halphen pencils
on quartic threefolds. Adv. Math. {\bf 223} (2010), no. 2,
594-618.\vspace{0.3cm}

\noindent 8. Cheltsov I.A., Nonexistence of elliptic structures
on general Fano complete intersections of index one. Moscow Univ.
Math. Bull.  {\bf 60} (2005), no. 3, 30-33.\vspace{0.3cm}

\noindent 9. Pukhlikov A.V., Birationally rigid Fano
hypersurfaces, Izvestiya: Mathematics. {\bf 66} (2002), no. 6,
1243-1269.\vspace{0.3cm}

\begin{flushleft}
\it pukh@liv.ac.uk \\
\it pukh@mi.ras.ru
\end{flushleft}

\end{document}